\documentclass[a4paper,10pt]{amsart}

\usepackage[latin1]{inputenc}

\usepackage{latexsym}

\usepackage{amssymb}
\usepackage{amsfonts}
\usepackage{amsmath}
\usepackage{amsthm}
\usepackage{textcomp}

\usepackage{comment}

\newtheorem{theorem}{Theorem}[section]    
\newtheorem{lemma}[theorem]{Lemma}           
                                     
\newtheorem{corollary}[theorem]{Corollary}
\newtheorem{proposition}[theorem]{Proposition}

\theoremstyle{definition}

\title[Lempert's elliptic tubes]{Convexity properties and complete hyperbolicity of Lempert's elliptic tubes}

\author[D. Alessandrini]{Daniele Alessandrini}
\address{Daniele Alessandrini\\Institut de Recherche Math\'ematique Avanc\'ee\\CNRS and UDS\\7 rue Ren\'e Descartes\\F-67084 Strasbourg Cedex\\France}
\email{daniele.alessandrini@gmail.com}

\author[A. Saracco]{Alberto Saracco}
\address{Alberto Saracco\\Dipartimento di Matematica\\Universit\`a di Parma\\Viale Usberti 53/A\\I-43100 Parma\\Italy}
\email{alberto.saracco@unipr.it}

\date{\today}


\newcommand{\nuovo}[1]{{{\bfseries \upshape #1}}}

\newcommand{\enne}{{\mathbb{N}}}

\newcommand{\erre}{{\mathbb{R}}}
\newcommand{\ci}{{\mathbb{C}}}
\newcommand{\cp}{{\mathbb{CP}}}
\newcommand{\pro}{{\mathbb{P}}}
\newcommand{\cappa}{{\mathbb{K}}}
\newcommand{\iper}{{\mathbb{H}}}
\newcommand{\famil}{{\mathcal{F}}}
\DeclareMathOperator{\ann}{Ann}
\DeclareMathOperator{\interior}{int}

\newcommand{\freccia}{{\ \longrightarrow \ }}

\newcommand{\tende}{{\ \rightarrow \ }}

\begin{document}

\sloppy

\begin{abstract}
We prove that elliptic tubes over properly convex domains $D\subset\erre\pro^n$ are $\ci$-convex and complete Kobayashi-hyperbolic. We also study a natural construction of complexification of convex real projective manifolds.
\end{abstract}

\maketitle

\section{Introduction}

The notion of elliptic tube of a subset of $\erre^n$ was defined and studied by Lempert in \cite{Le}. For the definition, see section \ref{sec:tubes}, where we deal with subsets of $\erre\pro^n$ to emphasize the projective invariance of the construction.

If $D \subset \erre\pro^n$, the elliptic tube $D^e$ is a subset of $\ci\pro^n$, a sort of ``complexification'' of $E$. Some properties of the set $D$ transfer automatically to $D^e$, for example openness, closedness, connectedness and boundedness of $D$ imply the same for $D^e$.

The same is not true for convexity, for example Lempert shows that the elliptic tube over a triangle in $\erre^2$ is not a convex subset of $\ci^2$. This is not too surprising, as the property of convexity in $\ci^n$ is not projectively invariant.  
In \cite{Le}, Lempert was particularly interested in the case where the base domain $D$ is a properly convex domain. He showed that in this case the elliptic tube $D^e$ is linearly convex, hence pseudoconvex. He also proved that the Hilbert distance on $D$ is the restriction of the Kobayashi distance on $D^e$. In \cite{TV} a Monge-Amp\`ere maximal problem is studied on the elliptic tube $D^e$ when $D$ is properly convex.

In this paper we show that if $D$ is a properly convex domain, the elliptic tube $D^e$ is $\ci$-convex, and that the Kobayashi distance on $D^e$ is complete. We use the projective invariance of the elliptic tubes to construct a natural complexification of convex real projective manifolds, a geometric structure in the sense of Klein.

Section \ref{sec:lincon} introduces the different notions of convexity needed in the sequel, namely linear convexity and $\ci$-convexity, and the necessary instruments for the proofs in the following sections, among which an important one is the notion of projective dual complement. 

In section \ref{sec:tubes} are given the definition and the first properties of elliptic tubes. Then we prove that if $D \subset \erre\pro^n$ is an open convex set, then the tube over the dual complement of $D$ coincide with the dual complement of the tube over $D$. Hence $D^e$ is linearly convex. As a corollary, elliptic tubes over convex domains are pseudoconvex, holomorphically convex, domains of holomorphy, polynomially convex, Runge domains, and are convex with respect to the linear fractions.

In section \ref{sec:regular} we show some regularity properties of elliptic tubes. Namely, if $\partial D$ is of class $\mathcal C^k$ then also $\partial D^e\setminus \partial D$ is. Unless $D$ is very special (projectively equivalent to the ball), there is no regularity at the real points of the boundary of the elliptic tube. This prevents to use the equivalence between linearly convexity and $\ci$-convexity which holds for connected $\mathcal C^1$-smooth domains.
Then we prove that elliptic tubes over properly convex domains are $\ci$-convex. Using a characterization of Kobayashi complete hyperbolicity of $\ci$-convex domains given in \cite{NS}, this shows that elliptic tubes over properly convex domains are complete Kobayashi-hyperbolic, and that they are taut and hyperconvex.

In section \ref{sec:hyptubes} we apply our results to real and complex projective manifolds (in the sense of Klein). To every convex real projective manifold we associate in a natural way a complex projective manifold, that we call its complexification. We show that this complex projective manifold is complete Kobayashi-hyperbolic and that it is homeomorphic in a natural way to the tangent bundle of the original convex real projective manifold.

\vspace{0.2cm}

\emph{Acknowledgements.} We wish to thank Giuseppe Tomassini, Stefano Trapani and Sergio Venturini for helpful discussions.

A. Saracco was partially supported by the MIUR project \lq\lq Geometric properties
of real and complex manifolds\rq\rq.

\section{Linearly convex and $\ci$-convex sets}\label{sec:lincon}

In this section, we recall some results about linear convexity and $\ci$-convexity that we will need in the following. We will also recall the definition of convexity for a subset of the real projective space.

Let $\cappa$ be a field, $V$ be a finite dimensional vector space over $\cappa$, $\pro = \pro(V)$ be the corresponding projective space, $V^*$ be the dual and $\pro^* = \pro(V^*)$ be the projective dual. 

A subset $E \subset \pro$ is said to be \nuovo{linearly convex} if the complement of $E$ is a union of projective hyperplanes or, in other words, if every point in the complement of $E$ is contained in a projective hyperplane disjoint from $E$. 

The intersection of linearly convex sets is linearly convex, and the space $\pro$ is linearly convex. If $E \subset \pro$ is any set, the \nuovo{linearly convex hull} of $E$ is the intersection of all the linearly convex sets containing $E$. This is just the complement of the union of all the projective hyperplanes disjoint from $E$.
Every element $\xi \in \pro^*$ is a projective class of linear functionals, hence it has a well defined \nuovo{kernel} denoted by $\ker(\xi)$, a projective hyperplane of $\pro$. If $x \in \pro$, the set 
$$\ann(x) = \{ \xi \in \pro^* \ |\ x \in \ker(\xi) \}$$
is called the \nuovo{annihilator} of $x$, and it is a projective hyperplane of $\pro^*$. If $\phi:\pro \mapsto \pro^{**}$ is the canonical identification, then $\ann(x) = \ker(\phi(x))$.

Let $E \subset \pro$ be a set. Then the \nuovo{dual complement} of $E$ is the set
$$E^* = \{ \xi \in \pro^* \ |\ \ker(\xi) \cap E = \emptyset \} $$
i.e. it is the set of all projective hyperplanes disjoint from $E$. Note that if $E_1 \subset E_2$, then $E_2^* \subset E_1^*$. 

Let $E \subsetneq \pro$ be an open connected set. If $a \in \partial E$, we denote by $\Gamma(a)$ the set of all tangent hyperplanes to $E$ at $a$. We have $\Gamma(a) = \ann(a) \cap E^*$.

\begin{lemma}
If $E \subset \pro$, then
$$E^* = \pro^* \setminus \bigcup_{x \in E} \ann(x)$$
In particular $E^*$ is linearly convex.
\end{lemma}
\begin{proof}
Let $\xi \in E^*$. If $x\in \pro$ is such that $\xi \in \ann(x)$, then $x \in \ker(\xi)$, hence $x \not\in E$. Therefore $\xi \in \pro^* \setminus \bigcup_{x \in E} \ann(x)$.  

Let $\xi \not\in E^*$, then there is $x \in E \cap \ker(\xi)$. Hence $\xi \in \ann(x)$, with $x \in E$.
\end{proof}

If $F \subset \pro^*$, we will identify $F^* \subset \pro^{**}$ with $\phi^{-1}(F^*) \subset \pro$, and we will write $F^* = \phi^{-1}(F^*) \subset \pro$. More explicitely we have
$$F^* = \{ x \in \pro \ |\ \ann(x) \cap F = \emptyset \} $$
By previous lemma, we have the alternative description
$$F^* = \pro \setminus \bigcup_{\xi \in F} \ker(\xi)$$   

\begin{lemma}
$E^{**}$ is the linearly convex hull of $E$. In particular $E \subset E^{**}$, and if $E$ is linearly convex, then $E = E^{**}$
\end{lemma}
\begin{proof}
$x \in E^{**}$ iff $\ann(x)$ is disjoint from $E^*$ iff every projective hyperplane containing $x$ intersects $E$. Hence the complement of $E^{**}$ is the union of the projective hyperplanes not intersecting $E$.  
\end{proof}

Let $\cappa = \erre$ or $\cappa = \ci$. It is easy to see that if $E$ is open, then $E^*$ is compact, and if $E$ is compact then $E^*$ is open. $E^*$ is bounded in some affine chart if and only if $E$ has interior points. 

Let $E \subset \pro$. A projective hyperplane is said to be \nuovo{tangent} to $E$ if it intersects the boundary $\partial E$ but it does not intersects the interior part $\interior(E)$.

\begin{proposition}  \label{prop:tangent}
Let $E \subset \pro$ be open. Then ${(\overline{E})}^* = \interior(E^*)$. Moreover we have
$$\forall \xi \in \pro^* : \xi \in \partial E^* \Leftrightarrow \ker(\xi) \mbox{ is tangent to } E $$
Let $E \subset \pro$ be compact. Then ${(\interior(E))}^* \supset \overline{E^*}$. Moreover we have
$$\forall \xi \in \pro^* : \xi \in \partial E^* \Rightarrow \ker(\xi) \mbox{ is tangent to } E  $$
\end{proposition}
\begin{proof}
See \cite[prop. 2.5.1]{PS}.
\end{proof}

If $\cappa = \erre$ a set $E \subset \pro$ is said to be \nuovo{convex} if it does not contain projective lines and if the intersection with every projective line is connected. An open or compact convex set is homeomorphic to an open or closed ball. 

\begin{theorem} The connected components of a linearly convex set are convex. Open or compact convex sets are linearly convex. If $E$ is connected and open or compact, then $E$ is convex iff it is linearly convex. Moreover, the dual complement of a convex set is convex. 
\end{theorem}
\begin{proof}
See \cite[prop. 1.3.4]{PS}, \cite[thm. 1.3.6]{PS} and \cite[thm. 1.3.11]{PS}
\end{proof}

A convex set $E \subset \erre\pro^n$ is said to be \nuovo{properly convex} if it is not contained in a projective hyperplane and it does not contain an affine line.

\begin{theorem} Let $E \subset \erre\pro^n$ be a properly convex open set. Then:
\begin{enumerate}
\item $E$ is relatively compact in an affine chart $\erre\pro^n \setminus H$.
\item $\interior(E)$ and $\overline{E}$ are also properly convex, and $\interior\left(\overline{E}\right) = \interior(E)$ and $\overline{\interior(E)} = \overline{E}$
\item ${(\overline{E})}^* = \interior\left(E^*\right)$ and ${(\interior(E))}^* = \overline{E^*}$.
\end{enumerate}
\end{theorem} 
\begin{proof}
See \cite[prop. 1.3.1]{PS} and \cite[thm. 1.3.14]{PS}.
\end{proof}

\begin{proposition}  \label{prop:convex open neighborhoods}
Let $E \subset \erre\pro^n$ be a compact convex set. Then $E$ has a basis of properly convex neighborhoods.
\end{proposition}
\begin{proof}
A linearly convex set has a basis of linearly convex neighborhoods (see \cite[page 17]{PS}) $(U_\alpha)$. As $E$ is connected, it is always contained in the interior of a connected component of $U_\alpha$, that is convex. 
\end{proof}

\begin{proposition} Let $D\subset\erre^n$ be a convex open set not containing straight lines. Then for all $\varepsilon>0$, there are two strictly convex domains $D^1=D^1(\varepsilon)$ and $D^2=D^2(\varepsilon)$ with real analytic boundary, such that $D^1\subset D\subset D^2$, and $D^2\setminus D^1$ is contained in the $\varepsilon$-neighbourhood of $\partial D$.
\end{proposition}
\begin{proof} It is theorem 2.1 in \cite{PT}.
\end{proof}

\begin{corollary}\label{esaustione} Let $D\subset\erre\pro^n$ be a properly convex open set. Then there is an exhaustion of $D$ made by strictly convex open sets $D_k$ with real analytic boundary.
\end{corollary}
\begin{proof} Since $D$ is properly convex and open, is starred with respect to all of its points. Consider an affine chart containing $D$ such that $D$ is convex and $0\in D$. Fix a sequence $1>\delta_k>0$ strictly decreasing to $0$. We define $E_{1-\delta_k}=(1-\delta_k) D\subset D$.

Suppose we have defined $D_k\subset E_{1-\delta_k}$, strictly convex and with real analytic boundary (and $D_j\Subset D_l$, for $j<l<\bar k$), for all $k< \bar k$. We need to construct $D_{\bar k}\subset E_{1-\delta_{\bar k}}$. Fix $\varepsilon_{\bar k}<\min\{\delta_{\bar k}, d(\partial D_{\bar k-1},\partial E_{1-\delta_{\bar k}})\}$. 
$E_{1-\delta_{\bar k}}$ is convex and bounded, so by the previous proposition we can find $D_{\bar k}=D_{\bar k}(\varepsilon_{\bar k})\subset E_{1-\delta_{\bar k}}$ strictly convex with real analytic boundary such that $E_{1-\delta_{\bar k}}\setminus D_{\bar k}$ is contained in the $\varepsilon_{\bar k}$-neighbourhood of $E_{1-\delta_{\bar k}}$. So $D_{\bar k}\Supset D_{\bar k -1}$.

The defined sets $D_k$ have the required regularity properties and form an exhaustion of $D$.
\end{proof}

If $\cappa=\ci$ a set $E \subset \pro$ is \nuovo{$\ci$-convex} if it does not contain projective lines and if the intersection with every projective line is connected and simply connected. An open $\ci$-convex set is homeomorphic to an open ball (see \cite[thm. 2.4.2]{PS}).

\begin{theorem} If $E \subsetneq \cp^n$ is an open connected set, with $n > 1$, $E$ is $\ci$-convex iff for all $a \in \partial E$, $\Gamma(a)$ is non-empty and connected. 
\end{theorem}
\begin{proof} See \cite[thm. 2.5.2]{PS}.
\end{proof}

\section{Elliptic tubes}\label{sec:tubes}

We denote by $\pi:\erre^{n+1} \setminus \{0\} \mapsto \erre\pro^n$ the natural projection.

Let $I \subset \erre\pro^1$ be an interval, with extremes $a_0$ and $a_1$. Then there are two vectors $v_0, v_1 \in \erre^2$ such that $\pi(v_0) = a_0$, $\pi(v_1) = a_1$ and 
$$\interior(I) = \pi(\{ c_0 v_0 + c_1 v_1 \ |\ c_0, c_1 \in \erre, c_0 c_1 > 0 \}) $$
$$\overline{I} = \pi(\{ c_0 v_0 + c_1 v_1 \ |\ c_0, c_1 \in \erre, c_0 c_1 \geq 0 \})$$
Consider the dual basis $v^0, v^1 \in {(\erre^2)}^*$. Now we have
$$\interior(I) = \{x \in \erre\pro^1 \ |\ v^1(x) v^0(x) > 0 \}$$
$$\overline{I} = \{x \in \erre\pro^1 \ |\ v^1(x) v^0(x) \geq 0 \}$$

Consider the natural inclusion $\erre\pro^1 \subset \cp^1$. Consider the sets
$${\interior(I)}^e = \pi(\{ c_0 v_0 + c_1 v_1 \ |\ c_0, c_1 \in \ci, \Re\left(c_0 \overline{c_1}\right) > 0 \}) = $$ $$= \{z \in \cp^1 \ |\ \Re\left(v^0(z) \overline{v^1(z)} \right) > 0 \} $$
$${\overline{I}}^e = \pi(\{ c_0 v_0 + c_1 v_1 \ |\ c_0, c_1 \in \ci, \Re\left(c_0 \overline{c_1}\right) \geq 0 \}) =$$ $$=  \{z \in \cp^1 \ |\ \Re\left(v^0(z) \overline{v^1(z)} \right) \geq 0 \} $$

It is easy to see that the sets ${\interior(I)}^e$ and ${\overline{I}}^e$ are just the open and closed circle in $\cp^1$ with diameter $I$. The circles ${\interior(I)}^e$ and ${\overline{I}}^e$ are well defined up to projective changes of coordinates in $\erre\pro^1$. 

In the same way, if $I \subset \erre\pro^n$ is a segment with extremes $a_0$ and $a_1$ there are two vectors $v_0, v_1 \in \erre^{n+1}$ such that $\pi(v_0) = a_0$, $\pi(v_1) = a_1$ and 
$$\interior(I) = \pi(\{ c_0 v_0 + c_1 v_1 \ |\ c_0, c_1 \in \erre, c_0 c_1 > 0 \}) $$
$$\overline{I} = \pi(\{ c_0 v_0 + c_1 v_1 \ |\ c_0, c_1 \in \erre, c_0 c_1 \geq 0 \})$$

If $L$ is the real projective line containing $I$, we denote by $L^\ci$ the complexification of $L$. Then ${\interior(I)}^e$ and ${\overline{I}}^e$ are the circles in $L^\ci$ with diameter $I$, well defined as in the previous section:
$${\interior(I)}^e = \pi(\{ c_0 v_0 + c_1 v_1 \ |\ c_0, c_1 \in \ci, \Re\left(c_0 \overline{c_1}\right) > 0 \})$$
$${\overline{I}}^e = \pi(\{ c_0 v_0 + c_1 v_1 \ |\ c_0, c_1 \in \ci, \Re\left(c_0 \overline{c_1}\right) \geq 0 \})$$

Let $D \subset \erre\pro^n$ be a set. The \nuovo{elliptic tube} with base $D$ is defined (as in \cite{Le}) as the set
$$D^e = \bigcup \{ I^e \ |\ I \subset D \mbox{ is a closed segment } \} \subset \cp^n $$
The definition is projectively invariant. If $D$ is open, then this definition is equivalent to
$$D^e = \bigcup \{ I^e \ |\ I \subset D \mbox{ is an open segment } \} \subset \cp^n  $$

Note that if $D$ is open also $D^e$ is open, and if $D$ is closed also $D^e$ is closed. If $D$ is connected also $D^e$ is connected. If $D$ is contained in an affine chart also $D^e$ is contained in the same affine chart, and if $D$ is bounded in an affine chart, also $D^e$ is bounded in the same affine chart. Moreover, if $D$ is contained in an affine chart, then $D$ is a deformation retract of $D^e$. If $D$ is open (non necessarily affine) there is a projectively invariant deformation retraction from $D^e$ to $D$. In these cases $D^e$ is homotopically equivalent to $D$.

Let $D \subset \erre\pro^n$ be an open convex set. Then $\pi^{-1}(D) \subset \erre^{n+1}$ is the union of two disjoint open convex cones. We choose one of these convex cones, and we denote it by $\tilde{D}$.

Let $D^* \subset {(\erre\pro^n)}^*$ be the dual complement of $D$. Then $\pi^{-1}(D) \subset {(\erre^{n+1})}^*$ is the union of two disjoint open convex cones. Exactly one of these two convex cones contains only linear functionals that are positive on $\tilde{D}$, we denote this cone by ${(\tilde{D})}^*$. 

Let $\famil \subset {(\erre\pro^n)}^*$ be a compact set such that $\famil^* = D$. In other words $\famil$ is a compact set such that the linearly convex hull of $\famil$ is $D^*$. If $\xi \in \famil$, we denote by $\tilde{\xi}$ an element of ${(\tilde{D})}^*$ such that $\pi(\tilde{\xi}) = \xi$. Then 
$$D = \{ x \in \erre\pro^n \ |\ \forall f,g \in \famil: \tilde{f}(x)\tilde{g}(x) > 0 \} $$

\begin{theorem}
Let $D$ and $\famil$ as above. Then
$$D^e = \{ z \in \cp^n \ |\ \forall f,g \in \famil: \Re\left(\tilde{f}(z)\overline{\tilde{g}(z)}\right) > 0 \} $$
\end{theorem}
\begin{proof}
See \cite[thm. 2.1]{Le}.
\end{proof}

\begin{theorem}
Let $D$ and $\famil$ as above. Then
$${(\famil^e)}^* = D^e$$
\end{theorem}
\begin{proof}
${(\famil^e)}^* \supset D^e$: We have to show that if $\xi \in \famil^e$, then $\ker(\xi) \cap D^e = \emptyset$. If $\xi \in \famil^e$, then by definition of elliptic tubes there are $f,g \in \famil$ such that $\tilde{\xi} = c_0 \tilde{f} + c_1 \tilde{g}$ with $c_0, c_1 \in \ci$, $\Re\left(c_0 \overline{c_1}\right) \geq 0$. Let $z\in \cp^n$ be such that $\tilde{\xi}(z) = c_0 \tilde{f}(z) + c_1 \tilde{g}(z) = 0$. Then 
$$\Re\left( \tilde{f}(z)\overline{\tilde{g}(z)} \right) = |\tilde{g}(z)| \Re\frac{\tilde{f}(z)}{\tilde{g}(z)} = - |\tilde{g}(z)| \Re\frac{c_1}{c_0} = - \frac{|\tilde{g}(z)|}{|c_0|} \Re\left(c_0 \overline{c_1}\right) \leq 0$$
Hence $z \not\in D^e$.

${(\famil^e)}^* \subset D^e$: We have to show that if $z \not\in D^e$ there exists an element $\xi \in \famil^e$ such that $z \in \ker(\xi)$. If $z \not\in D^e$ there are $f,g \in \famil$ such that $\Re\left(\tilde{f}(z)\overline{\tilde{g}(z)}\right) \leq 0$. Now if $\xi$ is such that $\tilde{\xi} = \tilde{g(z)} \tilde{f} - \tilde{f}(z) \tilde{g}$, then $\xi \in \famil^e$, and $z \in \ker(\xi)$.  
\end{proof}

\begin{corollary}\label{cor:lcopen}
If $D \subset \erre\pro^n$ is an open convex set, then $D^e$ is linearly convex. 
\end{corollary}
\begin{proof}
A dual complement is always linearly convex. Note that this statement also follows from the first part of the proof of \cite[thm. 2.2]{Le}. 
\end{proof}

\begin{lemma}
Let $D \subset \erre\pro^n$ be a compact convex set, and let $(U_k)$ be a family of open convex sets such that $U_{k+1} \subset U_k$ and $\bigcap U_k = D$. Then $D^e = \bigcap U_k^e$.
\end{lemma}
\begin{proof}
$D^e \subset \bigcap U_k^e$: For all $k$ we have $D \subset U_k$, hence $D^e \subset U_k^e$. 

$D^e \supset \bigcap U_k^e$: We only have to prove that if $L \subset \erre\pro^n$ is a real projective line, and $L^\ci$ is its complexification, then $D^e\cap L^\ci \supset L^\ci \cap \left(\bigcap U_k^e \right) = \bigcap \left( U_k^e \cap L^\ci \right)$. This is easy because $L \cap D$ and $L \cap U_k$ are just intervals, and $L^\ci \cap D^e$ and $L^\ci \cap U_k$ are just circles.  
\end{proof}

\begin{corollary}
If $D \subset \erre\pro^n$ is a compact convex set, then $D^e$ is linearly convex. 
\end{corollary}
\begin{proof}
By proposition \ref{prop:convex open neighborhoods}, we can always construct a family of open properly convex sets such that $U_{k+1} \subset U_k$ and $\bigcap U_k = D$. By previous lemma $D^e = \bigcap U_k^e$. If $z \not\in D^e$, we can find a $k$ such that $z \not\in U_k^e$, and as, by previous corollary, $U_k^e$ is linearly convex, then there is a projective hyperplane $H$ containing $z$ and disjoint from $U_k^e$. Hence $D^e$ is linearly convex.
\end{proof}

\begin{corollary} Let $D$ be an open or compact properly convex set. Then
$${(D^e)}^* = {(D^*)}^e$$
\end{corollary}
\begin{proof}
If $D$ is open, then by the previous theorem we know that $D^e = {({(D^*)}^e)}^*$. Hence ${(D^e)}^* = {({(D^*)}^e)}^{**} = {(D^*)}^e$, because by previous corollary ${(D^*)}^e$ is linearly convex.

If $D$ is compact, then $D^*$ is open, hence we have ${({(D^*)}^e)}^* = {(D^{**})}^e = D^e$. Hence ${(D^*)}^e = {(D^e)}^*$. 
\end{proof}

\begin{corollary}
Let $D$ be an open properly convex set. Then $D^e$ is pseudoconvex, holomorphically convex, a domain of holomorphy, polynomially convex, a Runge domain, and it is convex with respect to the linear fractions. Let $D$ be a compact properly convex set. Then $D^e$ is polynomially convex. 
\end{corollary}
\begin{proof}
This follows from the fact that $D^e$ is linearly convex and with connected dual complement. See \cite[prop. 2.1.8]{PS}, \cite[prop. 2.1.9]{PS}, \cite[prop. 2.1.11]{PS}.
\end{proof}

\section{Regularity and $\ci$-convexity of elliptic tubes}\label{sec:regular}

Let $D \subset \erre\pro^n$ be an open properly convex set, let $h_D$ be the Hilbert distance on $D$ (see \cite[sect. 3]{Le} for all the definitions needed here), and let $k_{D^e}$ be the Kobayashi distance on $D^e$. Note that as $D$ is bounded in some affine chart, also $D^e$ is bounded in the same affine chart, hence $D^e$ is Kobayashi hyperbolic, i.e. $k_{D^e}$ is a non degenerate distance.

\begin{theorem} If $D \subset \erre\pro^n$ is an open properly convex set, then 
\begin{enumerate}
\item For all $x,y \in D$, $h_D(x,y) = k_{D^e}(x,y)$.
\item If $L \subset \erre\pro^n$ is a real projective line, then $L^\ci \cap D^e$ is an open disk, and $k_{D^e}$ restricted to $L^\ci \cap D^e$ is the Poincar\'e distance on the disk.
\item If $z \in D^e$, let $L$ be the unique real line such that $L^\ci$ contains $z$ ($L^\ci$ is just the line containing $z$ and $\bar{z}$). Then there exists $x \in L$ such that 
$$k_{D^e}(z,x) = \min_{y \in D} k_{D^e}(z,y) $$
\end{enumerate}
\end{theorem}
\begin{proof}
Note that the first statement is \cite[thm. 3.1]{Le}. The other two can be proved in the same way. 
\end{proof}

Consider the functions:
$$d:D^e \ni z \mapsto \min_{y \in D} k_{D^e}(z,y) \in \erre_{\geq 0}$$
$$\phi:D^e \ni z \mapsto 2 \arctan(\tanh(d(z))) \in [0, \pi/2)$$

Choose a real projective hyperplane $H$ such that $D$ is bounded in $\erre^n = \erre\pro^n \setminus H$. Hence $D^e$ is bounded in  $\ci^n = \cp^n \setminus H^\ci$. With these coordinates, $D^e \subset D \oplus i \erre^n$. Consider the function
$$p: D \oplus i \erre^n \ni z=x+iy \mapsto \inf\{t > 0 \ |\ x+t^{-1}y \in D \} $$

\begin{lemma} Let $L\subset\erre^n$ be a real line. Then $L^\ci\cap (D\oplus i\erre^n)$ is a strip $S=(a,b) \oplus i\erre$ and 
$$p|_S(z) = p|_S (x+iy)= \begin{cases}\frac{y}{b-x} & {\rm if}\ y\geq0\\ \frac{y}{a-x} & {\rm if}\ y<0\end{cases}$$

Moreover,
$$\{p(z)p(\bar{z})=1\}=\partial D^e\setminus\partial D$$ 
\end{lemma}
\begin{proof} It is a simple computation.
\end{proof}

We can now define
$$u: D^e \ni z \mapsto \arctan\left(\frac{p(z)+p(\bar{z})}{1-p(z)p(\bar{z})}\right) $$

\begin{proposition}
$$u|_{D^e} = \phi$$
\end{proposition}
\begin{proof}
Choose a real line $L \subset \erre^n$, and compute explicitely $u$ on $L^\ci \cap D^e$. Also $\phi$ can be computed explicitely on $L^\ci \cap D^e$ using the previous theorem. Choose coordinates on $L^\ci \cap D^e$ such that it is the unit disk in $\ci$: then on this disk both functions are equal to
$$\arctan\left(\frac{2|\Im(z)|}{1-{|z|}^2}\right) = \left| \arg\frac{1+z}{1-z}\right| $$
\end{proof}

\begin{proposition}
Suppose that $D$ has a boundary of class $C^k$ (with $k \in \enne \cup \{\infty\} \cup \{\omega\}$ ). Then the functions $p$ and $u$ are of class $C^k$ on their domains. 
\end{proposition}
\begin{proof}
Let $z = x + iy \in D \oplus i \erre^n$. The half-line $\{x + t^{-1}y \ |\ t > 0\}$ cuts $\partial D$ in a single point $h \in \partial D$. By hypothesis, there exists a neighborhood $U$ of $h$ in $\ci^n$ and a function $f:U \mapsto \erre$ such that $\partial D \cap U = f^{-1}(0)$ and $\forall \zeta \in \partial D: df_\zeta \neq 0$. Then there exists a neighborhood $V$ of $z$ in $\ci^n$, such that for all $z' = x'+iy' \in V$, $p(w)$ is precisely the unique value of $t$ such that $f(x' + t^{-1}y') = 0$. Observe that
$$\frac{\partial}{\partial t} f(x+t^{-1}y) = -\frac1{t^2}\left< \nabla f(x+t^{-1}y),y\right>,$$
which on $\partial D$ is non vanishing since the direction $y$ is transversal to $\partial D$. Thus, by the implicit function theorem $p$ is of class $C^k$ near $z$.

Then also $u$ is of class $C^k$.
\end{proof}

\begin{proposition}
Suppose that $D$ has a boundary of class $C^k$ (with $k \in \enne \cup \{\infty\} \cup \{\omega\}$ ). Then $\partial D^e \setminus \partial D$ is of class $C^k$.
\end{proposition}
\begin{proof}
$\partial D^e \setminus \partial D \subset D \oplus i \erre^n$, and it is precisely the set where $p(z)p(\bar z) = 1$. Fix $\zeta\in \partial D^e\setminus\partial D$. Let $L^\ci$ be the complex line containing $\zeta$ and $\bar \zeta$. Then
$$p(z)p(\bar z)|_{L^\ci}=-\frac{y^2}{(b-x)(a-x)}$$
where $x+iy$ is the coordinate on $L^\ci$ corresponding to $z$. We have
$$\frac{\partial}{\partial y}(p(z)p(\bar z))|_{L^\ci\cap\partial D^e}=\frac2y\neq0$$ 
Hence by the implicit function theorem $\partial D^e\setminus \partial D$ is of class $C^k$.
\end{proof}

\begin{theorem}
Suppose that $D \subset \erre\pro^n$ is a properly convex open set with $\partial D$ of class $C^1$. Then $D^e$ is $\ci$-convex.
\end{theorem}
\begin{proof}
Let $a \in \partial D^e$. Since $D^e$ is linearly convex, $\Gamma(a)\neq\emptyset$. If $a \not\in \partial D$, we know that $\partial D^e$ is smooth at $a$, hence it has only one tangent real hyperplane $H$ at $a$. Note that $H$ has real dimension $2n-1$. It contains only one complex hyperplane, hence $\Gamma(a)$ is a single point, hence it is connected. If $a \in \partial D$, then $a$ is a real point. $\Gamma(a) = \ann(a) \cap {(D^e)}^* = \ann(a) \cap {(D^*)}^e$. By the second part of proposition \ref{prop:tangent}, $\ann(a) \cap \erre^n$ is a real hyperplane (of dimension $n-1$) tangent to $D$, hence $\ann(a)$ intersects ${(D^*)}^e$ in a single point, and also $\Gamma(a)$ is a single point hence it is connected.
\end{proof}

\begin{corollary}\label{bdCcvx}
Suppose that $D \subset \erre\pro^n$ is a properly convex open set. Then $D^e$ is $\ci$-convex.
\end{corollary}
\begin{proof}
By corollary \ref{esaustione}, there is an exhaustion of $D$ made by (strictly) convex domains $D_k$ whose boundary is real analytic. By theorem 5.1, each $D_k^e$ is $\ci$-convex. From the definition of the elliptic tubes is obvious that
$$D=\bigcup D_k\ \ \Rightarrow\ \ D^e=\bigcup D_k^e.$$
Hence $D^e$ is an increasing union of $\ci$-convex open sets, thus $\ci$-convex (see \cite[prop. 2.2.2]{PS}).
\end{proof}

\begin{corollary}\label{cor:tubecomphyp}
Suppose $D\subset\erre\pro^n$ is a properly convex open set. Then $D^e$ is complete Kobayashi hyperbolic, it is taut, and it is hyperconvex.
\end{corollary}
\begin{proof}
By previous corollary, $D^e$ is a bounded $\ci$-convex open set. Hence the statement follows from \cite[thm. 1]{NS}. 
\end{proof}

\begin{corollary}
Suppose that $D \subset \erre\pro^n$ is a convex open set. Then $D^e$ is $\ci$-convex.
\end{corollary}
\begin{proof}
If $D$ is bounded, it is corollary \ref{bdCcvx}. If $D$ is unbounded, it is increasing union of bounded convex open sets, hence $D^e$ is increasing union of $\ci$-convex open sets, hence $\ci$-convex.
\end{proof}

\begin{corollary}
Suppose that $D \subset \erre\pro^n$ is a compact convex set. Then $D^e$ is $\ci$-convex.
\end{corollary}
\begin{proof}
If $D$ is compact, then $D^*$ is open, hence ${(D^*)}^e$ is $\ci$-convex. By \cite[thm. 2.3.9]{PS} the dual complement of an open $\ci$-convex set is $\ci$-convex, hence $D^e$ is $\ci$-convex.
\end{proof}

\section{Complexification of convex real projective manifolds}\label{sec:hyptubes}

The projective invariance of the elliptic tubes may be used to construct a complexification of real manifolds with a suitable structure, namely with a structure of convex real projective manifold. The complexifications will have a structure of complex projective manifolds. These structures are geometric structures in the sense of Klein, see \cite{G} for a survey paper. Here we recall the definitions we need in the following. 

Let $\cappa$ be $\erre$ or $\ci$, and let $M$ be a manifold of dimension $n$ if $\cappa = \erre$, and of dimension $2n$ if $\cappa = \ci$. A $\cappa\pro^n$-structure on $M$ is given by a maximal atlas $\{(U_i,\phi_i)\}$, where the sets $U_i$ form an open covering of $M$ and the charts $\phi_i:U_i \freccia \cappa\pro^n$ are projectively compatible, i.e. the transition maps
$$\phi_{i,j} = \phi_i|_{U_i\cap U_j} \circ {\phi_j^{-1}}|_{\phi_j(U_i\cap U_j)} :\phi_j(U_i\cap U_j) \freccia \phi_i(U_i\cap U_j)  $$ 
have the property that for every connected component $C$ of the intersection $U_i \cap U_j$ there exists a projective map $A \in PGL_{n+1}(\cappa)$ such that $\phi_{i,j}|_C = A|_C$.
A $\cappa\pro^n$-manifold is a manifold with a $\cappa\pro^n$-structure. They are called real projective manifolds if $\cappa = \erre$ and complex projective manifolds if $\cappa = \ci$.

For example every open subset of $\cappa\pro^n$ (included $\cappa\pro^n$ itself) has a natural $\cappa\pro^n$-structure given by the inclusion map and all the charts that are compatible with the inclusion map. More interesting examples can be constructed by taking an open subset $\Omega$ of $\cappa\pro^n$ and a subgroup $\Gamma \subset PGL_{n+1}(\cappa)$ acting freely and properly discontinuously on $\Omega$. Then the quotient space $\Omega / \Gamma$ is a manifold and it inherits a $\cappa\pro^n$-structure from $\Omega$. The $\cappa\pro^n$-manifolds we will consider here are of this form.

The most interesting case is when $\cappa = \erre$ and $\Omega \subset \erre\pro^n$ is an open properly convex set. Real projective manifolds of the form $\Omega / \Gamma$ are called convex real projective manifolds. It is possible to construct many interesting manifolds of this form. For example, according to the Klein model of hyperbolic space, the hyperbolic space is identified with an ellipsoid $\iper^n \subset \erre\pro^n$, and the group of hyperbolic isometries is identified with the group of projective transformations of the ellipsoid, $O^+(1,n) \subset PGL_{n+1}(\erre)$. By this identification, every complete hyperbolic manifold is also a convex real projective manifold, and this gives plenty of interesting examples of convex real projective manifolds. It is also possible to construct many interesting examples when $\Omega$ is not an ellipsoid.  

Let $M = \Omega / \Gamma$ be a convex real projective manifold, in particular $\Omega \subset \erre\pro^n$ is an open properly convex set and $\Gamma \subset PGL_{n+1}(\erre)$ acts freely and properly discontinuously on $\Omega$. Consider the elliptic tube $\Omega^e$. By the projective invariance of the elliptic tube construction, the group $\Gamma$ also acts on $\Omega^e$.

\begin{proposition}
The action of $\Gamma$ on $\Omega^e$ is free and properly discontinuous.
\end{proposition}
\begin{proof}
To see that the action of $\Gamma$ on $\Omega^e$ is free, suppose, by contradiction, that an element $\gamma \in \Gamma$ has a fixed point $z \in \Omega^e \setminus \Omega$. As $\gamma$ is a real matrix, also the conjugate point $\bar{z}$ is fixed by $\gamma$, and also the unique complex line containing $z$ and $\bar{z}$. This complex line is real (in the sense of conjugation-invariant) hence it intersects $\Omega$ in a segment and it intersects $\Omega^e$ in a disc. As $\gamma$ does not fix any point of $\Omega$, it acts on the segment as a non-trivial translation. This action extends to the disc without fixing any point of the disc, and this is a contradiction with the fact that $z$ was a fixed point.  

Consider the group $G$ of all bi-holomorphisms of $\Omega^e$, equipped with the compact-open topology. As $\Omega^e$ is Kobayashi-hyperbolic, the group $G$ is a Lie group and it acts properly on $\Omega^e$ (see the introduction of \cite{Is} for a discussion of these properties). 

As $\Gamma$ acts properly discontinuously on $\Omega$, it is discrete for the topology it has as a subgroup of $PGL_{n+1}(\erre)$, that is the same it has as a subgroup of $PGL_{n+1}(\ci)$. The topology of this latter group is the compact-open topology for its action on $\ci\pro^n$ (see \cite[subsect. 2.6]{Go}). Hence the group $\Gamma$ is discrete even with the topology it has as a subgroup of $G$.  

This implies that $\Gamma$ is closed in $G$. In fact, by \cite[chap. 2, 1.8]{Pr}, a subgroup $H$ of $G$ is discrete if and only if there exists a neigborhood $U$ of $1$ in $G$ such that $H \cap U = \{1\}$. Now $\Gamma$ is discrete, its closure $\overline{\Gamma}$ is again a subgroup, but $\overline{\Gamma} \cap U = \{1\}$, hence $\overline{\Gamma}$ is discrete itself, hence $\Gamma = \overline{\Gamma}$.

As $\Gamma$ is closed in $G$ and the action of $G$ on $\Omega^e$ is proper, then also the action of $\Gamma$ on $\Omega^e$ is proper. As $\Gamma$ is discrete, the action is properly discontinuous.  
\end{proof}

We have seen that $\Gamma$ acts freely and properly discontinuously on $\Omega^e$. Hence the quotient $M^e = \Omega^e / \Gamma$ is a manifold and it has a natural complex projective structure. We call the complex projective manifold $M^e$ the \nuovo{complexification} of $M$. In the remaining part of this section we describe the manifold $M^e$. The inclusion $\Omega \subset \Omega^e$ gives an inclusion $M \subset M^e$. The complex conjugation on $\Omega^e$ is compatible with the action of $\Gamma$, hence it induces an anti-holomorphic involution on $M^e$ that has $M$ as locus of fixed points. The Kobayashi metric of $M^e$ is the quotient of the Kobayashi metric on $\Omega^e$ by the action of $\Gamma$, hence $M^e$ is a Kobayashi-hyperbolic complex manifold. As a corollary of the theorems in the first part we obtain

\begin{theorem}
Let $M$ be a convex real projective manifold. Then the Kobayashi metric on the complexification $M^e$ is complete. 
\end{theorem}
\begin{proof} It follows from the observations made above and corollary \ref{cor:tubecomphyp}.
\end{proof}

Finally, we describe the topology of $M^e$, by showing that it is homeomorphic in a natural way with the tangent bundle to $M$.

\begin{theorem}\label{Pepito}
Let $M$ be a convex real projective manifold, $M^e$ be its complexification and $TM$ be the tangent bundle of $M$. Then there is a natural homeomorphism
$$f: M^e \tilde{\rightarrow} TM$$
whose restriction to $M$ is the zero-section of $TM$.
\end{theorem}
\begin{proof} Let $M = \Omega / \Gamma$, with $\Omega \subset \erre\pro^n$ an open properly convex set and $\Gamma \subset PGL_{n+1}(\erre)$ acts freely and properly discontinuously on $\Omega$. We will give a natural homeomorphism  $f$ between $\Omega^e$ and the tangent space of $\Omega$, $T\Omega = \Omega \times \erre^n$, which is projectively invariant, hence passes to the quotient.

If $x\in \Omega$, then define $f(x)\in T\Omega$ as $(x,0)$. If $z\in \Omega^e\setminus \Omega$, consider the complex line $L_z$ through $z$ and $\bar{z}$. $I_z=L_z\cap \Omega$ is a segment, and $\Delta_z=L_z\cap \Omega^e$ is a disk with diameter $I_z$. Consider now the geodesic for the Poincar\'e metric $\gamma_z$ joining $z$ and $\overline z$ in $\Delta_z$. This is also a geodesic for the Kobayashi metric in $\Omega^e$. Define $x_z=\gamma_z \cap \Omega=\gamma_z\cap I_z$. We will define $f(z)$ as a vector in the tangent space at $x_z$, $T_{x_z}\Omega$. Since $\Omega^e$ is a complex manifold, it has a complex structure $J:T\Omega^e\to T\Omega^e$. Let us denote by $v\in T_{x_z}\Omega^e$ the unitary tangent vector to $\gamma$ at $x_z$ (considering $\gamma$ as a curve from $z$ to $\bar{z}$). Since $\gamma$ is a geodesic connecting two complex conjugate points, $Jv\in T_{x_z}\Omega$. Note that $Jv$ is tangent to $I_z$, thus the complex structure allows us to choose a direction on $I_z$ at $x_z$. Let $l_z=k_{\Omega^e}(z,x_z)$ (the Kobayashi distance). Then we can define
$$f(z)\ =\ (x_z,l_z\cdot Jv)\ \in\ T\Omega\, .$$

Note that the construction of $f$ is projectively invariant. We need to show that $f$ is a homeomorphism.

$f$ is {\em surjective}. Let $(x,w)\in T\Omega$. Consider the interval 
$$I_{(x,w)}=\Omega\cap\{x+tw\,|\,t\in\erre\}\, ,$$
and the elliptic tube over it, $\Delta_{(x,w)}=I^e_{(x,w)}$. Up to a projective trasformation we may suppose $\Delta_{(x,w)}=\Delta$, the unit disk, and $x=0$. Consider the two imaginary points $z_1, z_2$ such that $k_\Delta(x,z_i) = \Vert w\Vert$. Then we have that $f(z_1)=(x,\pm w)$ and $f(z_2)$ is the opposite point. 

$f$ is {\em injective}. Suppose $z_1,z_2\in \Omega^e$ are such that $f(z_1)=f(z_2)$. Then the complex geodesics $\gamma_{z_1}$ and $\gamma_{z_2}$ must coincide, they are at the same distance from the real part $\Omega$ and on the same side (in the disk, where the real part disconnects), hence $z_1=z_2$.

$f$ is {\em continuous}. Let us consider a sequence $\{z_n\}_{n\in\enne}\subset \Omega^e$ with $$\lim_{n\to\infty}z_n=z_\infty\in \Omega^e\,.$$ We have to show that $$\lim_{n\to\infty}f(z_n) = f(z_\infty)\,.$$ 

Note that the points $\overline z_n$ converge to $\overline z_\infty$ and the Kobayashi distances $k_{\Omega^e}(z_n,\overline z_n) $ converge to $k_{\Omega^e}(z_\infty,\overline z_\infty)$.

First consider the case where $z_\infty \in \Omega$. Consider the points $x_{z_n}=\gamma_{z_n} \cap \Omega$. In this case we just need to estimate the distance $k_{\Omega^e}(x_{z_n}, z_\infty) \leq k_{\Omega^e}(x_{z_n}, z_n) + k_{\Omega^e}(z_n, z_\infty) \tende 0$.

Then consider the case where $z_\infty \in \Omega^e \setminus \Omega$. In this case we can suppose that all the points $z_n$ are in $\Omega^e \setminus \Omega$, then the disks $\Delta_{z_n}$ are well defined and they converge to the disk $\Delta_{z_\infty}$. All we need to prove is that the part of the geodesics $\gamma_n$ connecting $z_n$ and $\overline z_n$ converge to the part of the geodesic $\gamma_\infty$ connecting $z_\infty$ and $\overline z_\infty$. Consider the limit set of the geodesics $\gamma_n$, defined as
$$\Lambda = \{t\in TM \ |\ t=\lim t_k \mbox{ for some sequence } t_k\in\gamma_{n_k}\}$$
Observe that $\Lambda$ is connected since it is the limit set of connected sets in a finite dimensional euclidean space, that $z_\infty,\overline z_\infty\in\Lambda$ since they are the limit of the points $z_k$ and $\overline z_k$ respectively, and that $\Lambda\subset\Delta_{z_\infty}$. We have to prove that that $\Lambda$ coincides with $\gamma_\infty$. 

We first prove that $\Lambda\subset\gamma_\infty$. Arguing by contradiction, let us suppose that there is a point $\zeta\in\Lambda\setminus\gamma_\infty$. By the unicity of the geodesic in the disk this means
\begin{eqnarray*}k_{\Omega^e}(z_\infty,\overline z_\infty) &=& k_{\Delta_\infty}(z_\infty,\overline z_\infty)<\\ &<& k_{\Delta_\infty}(z_\infty,\zeta)+k_{\Delta_\infty}(\zeta,\overline z_\infty)=k_{\Omega^e}(z_\infty,\zeta)+k_{\Omega^e}(\zeta,\overline z_\infty)\, .\end{eqnarray*} 
Since $\zeta\in\Lambda$, this means that 
$$\lim_{n\to\infty} k_{\Omega^e}(z_n,\overline z_n)\ \geq\ k_{\Omega^e}(z_\infty,\zeta)+k_{\Omega^e}(\zeta,\overline z_\infty)\ >\ k_{\Omega^e}(z_\infty,\overline z_\infty)\, ,$$
which is a contradiction. So $\Lambda\subset\gamma_\infty$.

Since $\Lambda$ is a connected set contained into $\gamma_\infty$ which is a topological closed interval and the endpoints of $\gamma_\infty$ belong to $\Lambda$, indeed $\Lambda=\gamma_\infty$.

We have proved that $f$ is a continuous bijective map between two domains of the same dimension. By Brouwer's invariance of domain theorem it is a homeomorphism.  
\end{proof}

\end{document}